\documentclass[reqno,tbtags,a4paper,12pt]{amsart}

\usepackage[T2A]{fontenc}     
\usepackage[utf8]{inputenc} 
\usepackage[in]{fullpage}
\usepackage[all]{xy}
\usepackage[unicode=true]{hyperref}

\usepackage{amsmath,amssymb,amsthm,amsfonts,amscd}

\def\colon{\mathpunct{:}}

\author{Vasilii Rozhdestvenskii}

\thanks{This work was supported by the Russian Science Foundation under grant no. 23-11-00143, \texttt{https:/\!/rscf.ru/project/23-11-00143/} }

\title{Domination of manifolds by hypersurfaces}
\date{}

\address{Steklov Mathematical Institute of Russian Academy of Sciences, Moscow, Russia}
\email{vrozhd@mi-ras.ru}

\begin{document}
\maketitle

\newtheorem{Th}{Theorem}
	\newtheorem*{Th*}{Theorem} 
	\newtheorem{Lem}{Lemma}
	\newtheorem*{Lem*}{Lemma}
	\newtheorem{example}{Example}
	\newtheorem*{example*}{Example}
	\newtheorem{Cor}{Corollary}
	\newtheorem{Cond}{Condition}
	\newtheorem*{Cor*}{Corollary}
	\newtheorem*{Frame}{Frame}
	\newtheorem*{Denom}{Denomination}
	\newtheorem*{Prop*}{Proposition}
	\newtheorem{Prop}{Proposition}
	\theoremstyle{definition}
	\newtheorem{definition}{Definition}
	\newtheorem*{definition*}{Definition}
	\newtheorem{Rem}{Remark}
	\newtheorem{Quest}{Question}
	\newtheorem{Conv}{Convention}

	\makeatletter
\newcommand{\colim@}[2]{%
  \vtop{\m@th\ialign{##\cr
    \hfil$#1\operator@font colim$\hfil\cr
    \noalign{\nointerlineskip\kern1.5\ex@}#2\cr
    \noalign{\nointerlineskip\kern-\ex@}\cr}}%
}
\newcommand{\colim}{%
  \mathop{\mathpalette\colim@{\rightarrowfill@\textstyle}}\nmlimits@
}
\makeatother
	
	\def\rk{\mathop{\mathrm{rk}}}
	\def\id{\mathord{\mathrm{id}}}
	\def\pt{\mathord{\mathrm{pt}}}
	\def\ker{\mathop{\mathrm{Ker}}}
	\def\coker{\mathop{\mathrm{Coker}}}
	\def\Im{\mathop{\mathrm{Im}}}
	\def\Tor{\mathop{\mathrm{Tor}}}
	\def\Ext{\mathop{\mathrm{Ext}}}
	\def\Tors{\mathop{\mathrm{Tors}}}
	\def\ord{\mathop{\mathrm{ord}}}
	\def\ex{\mathop{\mathrm{ex}}}
	\def\Indet{\mathop{\mathrm{Indet}}}
	\def\Hom{\mathop{\mathrm{Hom}}}
	\def\Ann{\mathop{\mathrm{Ann}}}
	\def\MU{\mathord{\mathrm{MU}}}
	\def\bu{\mathord{\mathrm{bu}}}
	\def\MSO{\mathord{\mathrm{MSO}}}
	\def\SO{\mathord{\mathrm{SO}}}
	\def\Sq{\mathord{\mathrm{Sq}}}
	\def\U{\mathord{\mathrm{U}}}
	\def\MG{\mathord{\mathrm{MG}}}
	\def\ch{\mathord{\mathrm{ch}}}
	\def\BU{\mathord{\mathrm{BU}}}
	\def\St{\mathord{\mathrm{St}}}
	\def\K{\mathord{\mathrm{K}}}
	\def\exp{\mathord{\mathrm{exp}}}
	\newcommand{\Z}{\mathbb{Z}}
	\newcommand{\Q}{\mathbb{Q}}
	\newcommand{\s}{\mathfrak{s}}
	\let\le\leqslant
	\let\ge\geqslant
	
	In this note, all manifolds are assumed to be $C^{\infty}$-smooth, connected, closed and orien\-ted, unless otherwise specified. Consider the set $\mathcal{M}_n$ of diffeomorphism types of mani\-folds of dimension $n$. This set can be endowed with a preorder $\gtrsim$ defined as follows: $M\gtrsim N$ if there exists a map $f\colon M \to N$ of  nonzero degree. Recall that maps of nonzero degree are called \textit{dominant} and the preorder is usually called the \textit{dominance relation} (for a general survey on the subject and further references, see \cite{Harpe}). For $n\le 2$ the dominance relation is completely understood: the set $\mathcal{M}_1$ consists of a single element and so the preorder is trivial; on $\mathcal{M}_2$ the preorder coincides with the ordering by the genus. The dominance relation on $\mathcal{M}_{3}$ is far more complicated, though many results about its structure are known. For $n\ge 4$ even general properties of the dominance relation on $\mathcal{M}_n$ are poorly understood.
	
	For any $n$, the sphere $S^n$ is a minimum in the proset $\mathcal{M}_n$. On the other hand, since dominant maps induce surjections on rational homology groups, for $n\ge2$, the preordered set $\mathcal{M}_n$ has no maximum. This motivates the following question. 

\begin{Quest}[Carlson---Toledo, 1989]\label{Q}
Is there an easily describable subset $\mathcal{C}_n\subset \mathcal{M}_n$ such that for any $M\in \mathcal{M}_n$ there exists $M'\in \mathcal{C}_n$ with $M'\gtrsim M$? 
\end{Quest}	

Of course, the question is posed rather informally and may admit different answers. First, let us mention two results that follow from the development of asphericalisation and hyperbolisation procedures, introduced by Kan and Thurston and by Gromov, respec\-tively: (1) \textit{one can take $\mathcal{C}_n$ to consist of all aspherical manifolds} (see \cite{DaJa}) and, more restrictively, (2) \textit{to consist of all manifolds which admit a Riemannian metric of negative sectional curvature in an arbitrarily small interval $[-1-\varepsilon, -1]$} (see \cite{Onta}).  Also, Kotschick and L\"oh conjectured that it would suffice \textit{to take $\mathcal{C}_n$ to consist only of hyperbolic manifolds}, i.e. those admitting a Riemannian metric of constant sectional curvature~$-1$. This conjecture is trivially true for $n=2$, was known for $n=3$ due to Brooks~\cite{Brooks}, was proved for $n=4$ by Gaifullin~\cite{Gaif2}, and remains open for $n\ge 5$.

 Another answer to Question~\ref{Q} was given by Gaifullin. From his explicit construction for realisation of cycles (see \cite{Gaif1}), it follows that \textit{one can take $\mathcal{C}_n$ to consist of all finite-sheeted coverings of the $n$-dimensional Tomei manifold}, i.e. the isospectral manifold of real symmetric tridiagonal $(n+1)\times (n+1)$ matrices. This led him to introduce the concept of a URC-manifold: a manifold $M$ is called a \textit{URC-manifold} if one can take $\mathcal{C}_n$ to consist of all finite-sheeted coverings of the manifold~$M$. In \cite{Gaif2} many more examples of URC-manifolds were constructed. 
 
The main theorem of this note provides one more possible answer to Question~\ref{Q}.  
 
\begin{Th}\label{Th}
Let $M$ be a manifold of dimension $n$. Then there exists a smooth, compact, connected submanifold $H\subset S^{n+1}$ of codimension $1$ and without boundary with a dominant map $f\colon H \to M$. Therefore, one can take $\mathcal{C}_n$ to consist of all codimension $1$ submanifolds of the $(n+1)$-dimensional sphere. 
\end{Th}
\begin{proof}
Recall that for any finite, connected $CW$-complex $X$, the suspension $\Sigma X$ is \linebreak rationally homotopy equivalent to a wedge sum of spheres (see \cite[Th. 2.2]{Berstein}).  A standard direct consequence of this equivalence is that the rational Hurewicz homomorphism $\mathrm{Hur}\colon \pi_*(\Sigma X)\otimes \Q \to H_*(\Sigma X; \Q)$ is surjective. Applying this to $X=M$, we obtain a map $g\colon S^{n+1} \to \Sigma M$ such that $g_*[S^{n+1}]=d\Sigma[M]$, where $d$ is a nonzero integer and square brackets denote the fundamental class. Perturbing $g$, if necessary, one can assume that~$g$ is transverse  to $M\subset \Sigma M$. Now put $H=g^{-1}(M)$ and $f=\left.g\right|_H$. Since the degree of a smooth map can be computed as the number of preimages (with signs) of a regular value, it follows that $d=\deg(f)=\deg(g)$. In order to make $H$ connected, one can join its connected components with tubes. 
\end{proof}

The following corollary is straightforward. 

\begin{Cor}
In every dimension $n$, there exists a codimension $1$ submanifold $H\subset S^{n+1}$ such that $H$ is a URC-manifold. 
\end{Cor}	
	
Finally, let us discuss a connection between Question \ref{Q} and the Steenrod problem on realisation of cycles. Recall that the Steenrod problem asks whether a given singular homology class $x\in H_n(X; \Z)$ of a topological space $X$ can be realised by a continuous image of the fundamental class of a manifold. From a famous result of Thom \cite[Th. III.4]{Thom} we know that any homology class can be realised in
the above sense after multiplying it by some nonzero integer. Therefore, if a subset $\mathcal{C}_n\subset \mathcal{M}_n$ serves as an answer to Question~\ref{Q}, then any $n$-dimensional integral homology class of a connected topological space can be realised with multiplicity by a manifold from $\mathcal{C}_n$. So, one gets the following corollary of Theorem~\ref{Th}. 

\begin{Cor}
Let $X$ be a connected topological space and $x\in H_n(X; \Z)$ be a homology class. Then there exists a codimension $1$ submanifold $H\subset S^{n+1}$ and a continuous map $f\colon H \to X$ such that $f_{*}[H]=dx$, where $d$ is a nonzero integer. 
\end{Cor}

\end{document}